\documentclass[format=acmsmall]{acmart}

\usepackage{amsmath}
\usepackage{amssymb}
\usepackage{amsfonts}
\usepackage{graphicx}
\usepackage{amsthm}
\usepackage{amscd}
\usepackage{amsmath}
\usepackage{latexsym}
\usepackage{xparse}
\usepackage{xifthen}
\usepackage{subfigure}
\usepackage[scaled=.8]{beramono}
\usepackage[boxed, lined]{algorithm2e}

\acmJournal{TOMS}
\acmVolume{V}
\acmNumber{N}
\acmArticle{A}
\acmYear{YYYY}
\acmMonth{1}
\copyrightyear{YYYY}

\setcopyright{acmcopyright}
\acmJournal{TOMS}
\acmYear{2019} \acmVolume{1} \acmNumber{1} \acmArticle{1} \acmMonth{1} \acmPrice{15.00}\acmDOI{10.1145/3321514}

\SetKwFor{ForM}{for}{}{end}
\SetKwProg{Fn}{function}{}{end}
\SetKwInOut{In}{input}
\SetKwInOut{Out}{output}

\newcommand{\Matlab}{\textsc{Matlab} }
\newcommand{\Matlabs}{\textsc{Matlab}}
\newcommand{\supp}{\mathrm{supp}}
\newcommand{\diag}{\mathrm{diag}}
\newcommand{\RR}{{\mathbb{R}}}

\newcommand{\NN}{{\mathbb{N}}}

\newcommand{\BSP}[3][]{{\mathcal{B}_{#1}[\ifthenelse{\isempty{#2}}{#3}{\seg{#2}{#3}}]}}
\DeclareDocumentCommand{\SPL}{ O{} O{} m m m}{{{\mathcal{#3}^{#2}[\mbf{#5},\mbf{#4}\ifthenelse{\isempty{#1}}{}{,\mbf{#1}}]}}}

\newcommand{\mbf}[1]{{\boldsymbol{#1}}}
\newcommand{\seg}[2]{{#2^{(#1)}}}
\newcommand{\per}[1]{{#1}_{\mathrm{per}}}
\newcommand{\bs}{b}
\newcommand{\bsext}{b}
\newcommand{\mdbs}{B}
\newcommand{\nseg}{{n_{s}}}
\newcommand{\ndof}{{n}}
\newcommand{\nsum}{{\mu}}
\newcommand{\mmult}{{m}}
\newcommand{\pdeg}{{p}}
\newcommand{\contin}{\kappa}
\newcommand{\deriv}{k}

\newcommand{\intI}{I}
\newcommand{\intIx}{x}

\newtheorem{proposition}{Proposition}
\newtheorem{remark}{Remark}
\newtheorem{example}{Example}
\newtheorem{definition}{Definition}

\begin{document}
\title{Algorithm xxx: Computation of Multi-Degree B-Splines}  

\author{Hendrik Speleers}
\affiliation{%
  \institution{Department of Mathematics, University of Rome Tor Vergata}
  \streetaddress{Via della Ricerca Scientifica~1}
  \city{Rome}
  \postcode{00133}
  \country{Italy}}
\email{speleers@mat.uniroma2.it}

\begin{abstract}
Multi-degree splines are smooth piecewise-polynomial functions where the pieces can have different degrees.
We describe a simple algorithmic construction of a set of basis functions for the space of multi-degree splines, with similar properties to standard B-splines. These basis functions are called multi-degree B-splines (or MDB-splines). The construction relies on an extraction operator that represents all MDB-splines as linear combinations of local B-splines of different degrees. 
This enables the use of existing efficient algorithms for B-spline evaluations and refinements in the context of multi-degree splines.
A \Matlab implementation is provided to illustrate the computation and use of MDB-splines.
\end{abstract}

%
%

\ccsdesc[500]{Mathematics of computing~Numerical analysis}
\ccsdesc[300]{Mathematics of computing~Computations on polynomials}%
%

\keywords{Multi-degree splines, MDB-splines, B-splines, Extraction operator}

\maketitle


\section{Introduction}

A \emph{spline} is a piecewise-polynomial function where the pieces are glued together in a smooth way. All the pieces are assumed to have the same polynomial degree. A \emph{B-spline} (or basis spline) is a spline that has minimal support with respect to a given degree, smoothness, and domain segmentation. Any spline of given degree can be represented as a linear combination of B-splines of that degree.
Thanks to their stable calculation and geometrically intuitive manipulation, B-spline representations form one of the foundations of the field of computer-aided geometric design (CAGD); see, for example, \cite{cohen2001}.
They are also a requisite ingredient in isogeometric analysis (IGA), a paradigm for the numerical
solution of partial differential equations; see, for example, \cite{cottrell2009}.

Generally, Schoenberg is regarded as the father of B-splines because of his seminal work on the topic \cite{schoenberg1946,currys1947}. It should be mentioned, however, that B-splines were already studied earlier (see, for example, the historical notes in \cite{butzer1988} and \cite{deboor2003}). Unfortunately, the B-spline construction based on divided differences of truncated power functions suffers from numerical instabilities.
A stable recurrence relation for the evaluation of B-splines was developed independently by \citet{deboor1972} and \citet{cox1972}. This recurrence relation opened the door for the use of B-splines for practical applications. The first software package for the stable calculation with B-splines was developed by  \citet{deboor1977} and many others followed afterwards. 
Nowadays, B-splines form an integral part of application fields like CAGD and IGA.

An extension of the standard concept of (B-)splines is to consider polynomial pieces of different degrees glued together with certain smoothness. This idea of \emph{multi-degree (B-)splines} dates back to \citet{nurnberger1984} but has been rejuvenated  several times under different forms (and names). 
A special case was investigated by \citet{sederberg2003knot} for degrees $(1,2,3)$ and $(1,n)$ in the context of CAGD. Later, \citet{shen2010changeable,shen2010md} constructed multi-degree B-splines with general degrees by means of an integral recurrence relation. \citet{li2012} observed the latter as being difficult to compute, and instead presented a geometric-oriented recursive procedure. 
Recently, multi-degree B-splines were revisited by \citet{beccari2017} illustrating their capabilities 
in CAGD, but still based on an integral recurrence relation and so-called transition functions identified through Hermite interpolation.
Even though the theoretical framework is very elegant, multi-degree B-splines are a burden to implement when dealing with recursive integration or when solving several Hermite interpolation problems. This might explain why they have not been popular in a practical context, despite their (theoretical) potential. 
Besides CAGD, multi-degree splines were also applied in IGA, in the frame of T-splines \cite{liu2016hybrid} and polar splines \cite{toshniwal2017polar}.
Non-polynomial extensions can be found in \cite{nurnberger1984,buchwald2003,sommer1988} considering local spaces generated by ECT-systems and Descartes systems.

The purpose of this paper is to show that the implementation of multi-degree B-splines (in short, MDB-splines) can be painless, and we provide a small \Matlab toolbox as illustration.
We closely follow the algorithmic construction proposed by \citet{toshniwal2017polar} in the context of rational multi-degree polar splines, but reformulate the algorithms in a computationally more efficient way. This further simplifies the construction, and might open the door to the practical use of multi-degree (B-)splines in a wide range of application fields, from CAGD and IGA to data fitting and compression. We remark that some computational aspects of multi-degree splines 
were also addressed in \cite{beccari2017}.

The remainder of the paper is organized as follows. We start by defining notation for standard B-splines in Section~\ref{sec:bsplines}. Section~\ref{sec:mdbsplines} introduces multi-degree splines as being composed of standard spline pieces that are joined by specified smoothness, and presents the algorithmic construction of a set of spanning MDB-splines with similar properties to B-splines. Section~\ref{sec:matlab} illustrates the small \Matlab toolbox with some examples, and we end in Section~\ref{sec:conclusion} with some concluding remarks.

\section{Uniform degree B-splines} \label{sec:bsplines}

In this section we define notation for standard B-splines (of uniform degree), and recall some of their main properties. We refer the reader to the monographs \cite{deboor2001,piegl2012nurbs} for more details.

Given a basic interval ${\intI} := [{\intIx_1},\;{\intIx_2}] \subset \RR$, let us denote with $U$ an \emph{open knot vector} of degree $\pdeg\in\NN$ and length $\ndof+\pdeg+1\in\NN$ (where $n\geq p+1$). 
More precisely, 
\begin{equation}
\begin{gathered}
  U := [{u_1},\;{u_2},\;\ldots,\;{u_{{\ndof+\pdeg+1}}}]\;,\quad {u_{i}}\leq {u_{i+1}},\\
  u_1 = \cdots = u_{{\pdeg}+1} = {\intIx_1} < u_{{\pdeg}+2},\quad
  u_{\ndof} < {\intIx_2} = u_{{\ndof+1}} = \cdots = u_{{\ndof+\pdeg+1}}.
\end{gathered}
\label{eq:knotVector}
\end{equation}
The number of times a knot value $u_i$ is duplicated in the knot vector is called the knot's \emph{multiplicity}. The multiplicity of $u_i$ is denoted by $\mmult_i$, and we assume that $1\leq \mmult_i \leq \pdeg+1$.
Such knot vector allows us to define $\ndof$ linearly independent B-splines of degree $\pdeg$.

\begin{definition}[B-splines]\label{def:bspline}
For a given open knot vector $U$ of degree $\pdeg$, the corresponding set of B-splines $\{\bs_{j,\pdeg}:j=1,\ldots,\ndof\}$ are defined using the recursive relation,
\begin{equation*}
  \bs_{j,\pdeg}(\intIx) := \frac{\intIx-u_j}{u_{j+\pdeg}-u_j}\bs_{j,\pdeg-1}(\intIx) + \frac{u_{j+\pdeg+1}-\intIx}{u_{j+\pdeg+1}-u_{j+1}}\bs_{j+1,\pdeg-1}(\intIx), \quad x\in\intI,
\end{equation*}
starting from
\begin{equation*}
  \bs_{j,0}(\intIx) := 
  \begin{cases}
    \; 1, & \text{if }u_j \leq \intIx < u_{j+1},\\
    \; 0, & \text{otherwise},
  \end{cases}
\end{equation*}
and under the convention that fractions with zero denominator have value zero.
\end{definition}

With the above definition, all the B-splines take the value zero at the end point $\intIx_2$. 
Therefore, in order to avoid asymmetry over the interval $\intI$, it is common to assume the B-splines to be left continuous at $\intIx_2$. We will follow suit.

\begin{proposition}\label{pro:properties-bspline}
The B-splines $\bs_{j,\pdeg}$, $j=1,\ldots,\ndof$ enjoy the following properties:
\begin{itemize}
  \item local support:
  \begin{equation*}
  \supp(\bs_{j,\pdeg})=[u_j,u_{j+\pdeg+1}], \quad j=1,\ldots,\ndof;
  \end{equation*}
  \item smoothness at a knot:
    \begin{equation*}
      \bs_{j,\pdeg} \in C^{\pdeg-\mmult_i}(u_i), \quad j=1,\ldots,\ndof;
    \end{equation*}
  \item non-negative partition of unity:
  \begin{equation*}
    \bs_{j,\pdeg}(\intIx)\geq0, \quad j=1,\ldots,\ndof, \quad \sum_{j=1}^{\ndof}\bs_{j,\pdeg}(\intIx)=1; 
  \end{equation*}
  \item interpolation at the end points:
  \begin{align*} 
   &\bs_{1,\pdeg}(\intIx_1)=1, \quad \bs_{j,\pdeg}(\intIx_1)=0, \quad j=2,\ldots,\ndof, \\
   &\bs_{\ndof,\pdeg}(\intIx_2)=1, \quad \bs_{j,\pdeg}(\intIx_2)=0, \quad j=1,\ldots,\ndof-1.
  \end{align*}
\end{itemize}
\end{proposition}

The spline space corresponding to $U$ is denoted with $\BSP{}{U}$ and is defined as the span of {$\{\bs_{j,\pdeg}:j=1,\dots,\ndof\}$}. This is a space of piecewise polynomials of degree $\pdeg$ with smoothness $C^{\pdeg-\mmult_i}$ at knot $u_i$ and its dimension is $\ndof$. 

We identify a spline function $f\in \BSP{}{U}$ with the vector of its coefficients $\left[f_1,\dots,f_\ndof\right]$, i.e.,
\begin{equation} \label{eq:splineRepr}
  f(\intIx) = \sum_{j=1}^{\ndof} f_j \bs_{j,\pdeg}(\intIx), \quad \intIx\in\intI.
\end{equation}
The $\deriv$-th derivative of $f$ (for $0\leq\deriv \leq \pdeg$) is given by
\begin{equation*}
D_+^\deriv f(\intIx) = \sum_{j=1}^{\ndof-\deriv} f_{j,\deriv} \bs_{j,\pdeg-\deriv}(\intIx), \quad \intIx\in\intI,
\end{equation*}
where the lower degree B-splines {$\bs_{j,\pdeg-\deriv}$} are defined over the smaller knot vector \linebreak
$[u_{\deriv+1},\;u_{\deriv+2},\;\dots,\;u_{n+p-\deriv+1}]$, and
\begin{equation*}
  f_{j,\deriv} := 
  \begin{cases}
  \; f_j, & \text{if }\deriv = 0,\\
  \; \alpha_{j,\deriv}(f_{j+1,\deriv-1}-f_{j,\deriv-1}), & \text{if }\deriv>0,
  \end{cases}
\end{equation*}
and 
\begin{equation*}
 \alpha_{j,\deriv} := \frac{\pdeg-\deriv+1}{u_{j+\pdeg+1}-u_{j+\deriv}}. 
\end{equation*}
In view of the interpolation property of B-splines at the end points of $\intI$, this implies that
only the first (last) $\deriv+1$ B-splines contribute towards the $\deriv$-th derivative at the left (right) end point of $\intI$. 
Hence, given $[f_1,\dots,f_\ndof]$, we can write the $\deriv$-th derivative of $f$ at the left or right end point of $\intI$ in terms of $[f_1,\dots,f_{\deriv+1}]$ or $[f_{\ndof-\deriv},\dots,f_{\ndof}]$, respectively.

\begin{remark} \label{rmk:structureU}
The structure of $U$ in \eqref{eq:knotVector} is such that $\pdeg$, $\mmult_i$, $\ndof$ and $\intI$ are embedded in it. Therefore, we can assume that once a knot vector is known, so are the degree, smoothness, and dimension of the corresponding spline space $\BSP{}{U}$.
\end{remark}

\begin{remark} \label{rmk:bezier}
When taking $U$ in \eqref{eq:knotVector} with $\ndof=\pdeg+1$, there are no knots in the interior of $\intI$ and the space $\BSP{}{U}$ coincides with the polynomial space of degree $\pdeg$. In this case, the B-spline basis is nothing other than the Bernstein polynomial basis, and the representation in \eqref{eq:splineRepr} is referred to as B\'ezier representation.
\end{remark}

\section{Multi-degree B-splines}\label{sec:mdbsplines}

In this section we focus on the multi-degree spline setting. In general, a multi-degree spline is a smooth piecewise-polynomial function where the pieces can have different degrees. 
First, we define our multi-degree spline space as a collection of standard spline spaces (with possibly different degrees) glued together smoothly. Then, we present an algorithmic construction of a set of basis functions for the multi-degree spline space, with similar properties to standard B-splines. These basis functions are called \emph{multi-degree B-splines} (in short, \emph{MDB-splines}). We closely follow the notation and presentation from \cite{toshniwal2017polar}, but in slightly simplified form. 
We detail the entire construction to keep the paper self-contained. 

\subsection{Multi-degree splines}

Consider $\nseg$ open knot vectors $\seg{i}{{U}}$ of degree $\seg{i}{\pdeg}$, $i=1,\dots,\nseg$, defined as in \eqref{eq:knotVector}. We denote the left and right end points of the interval $\seg{i}{\intI}$ associated to $\seg{i}{U}$ with $\seg{i}{\intIx_{1}}$ and $\seg{i}{\intIx_{2}}$, respectively. 
The collection $\mbf{U} := (\seg{1}{U},\dots,\seg{\nseg}{U})$ is called an \emph{$\nseg$-segment knot vector configuration}. 
The multi-degree spline spaces will be constructed by considering spline spaces over the knot vectors $\seg{i}{U}$, which are glued together with certain smoothness requirements at the end points $\seg{i}{\intIx_2}$ and $\seg{i+1}{\intIx_1}$ for $i \in \{1,2,\dots,\nseg-1\}$. The equivalence class at the points $\seg{i}{\intIx_2}$ and $\seg{i+1}{\intIx_1}$ is called the \emph{$i$-th segment join}. 
In the following $\NN_{-1}:= \NN\cup\{-1,0\}$.

\begin{definition}[Multi-degree spline space]\label{def:multiDegSpline}
  Given an $\nseg$-segment knot vector configuration $\mbf{U}$, we define for each segment the mapping
  \begin{equation*}
    \seg{i}{\phi}(\intIx) := \intIx-\seg{i}{\intIx_1} + \seg{1}{\xi_1} + \sum_{\ell=1}^{i-1}(\seg{\ell}{\intIx_2}-\seg{\ell}{\intIx_1}),
    \quad i=1,\ldots,\nseg.
  \end{equation*}
  Then, $\seg{i}{\Omega} := [\seg{i}{\xi_1},\;\seg{i}{\xi_2}] := \seg{i}{\phi}( [\seg{i}{\intIx_1},\;\seg{i}{\intIx_2}]) \subset \RR$, and we construct the composed interval
  \begin{equation*}
    \Omega := \seg{1}{\Omega}\cup\dots\cup\seg{\nseg}{\Omega}.
  \end{equation*}
  Moreover, let $\mbf{\contin} := (\seg{1}{\contin},\seg{2}{\contin},\ldots,\seg{\nseg-1}{\contin}) \in \NN_{-1}^{\nseg-1}$ be a given vector of continuity orders. The space of multi-degree splines is defined as
  \begin{equation*}
  \begin{split}
    \SPL{S}{\contin}{U} := \bigl\{\,&s:\Omega\rightarrow \RR: \;
    s\circ \seg{i}{\phi} \in \BSP{i}{U},\ \forall i\in\{1,\ldots,\nseg\} \;\;\& \\
    &D_-^\deriv( s\circ\seg{i}{\phi})(\seg{i}{\intIx_2}) 
    = D_+^\deriv (s\circ\seg{i+1}{\phi})(\seg{i+1}{\intIx_1}),\
      \forall \deriv \in \{0,\dots,\seg{i}{\contin}\},\ \forall i\in\{1,\ldots,\nseg-1\}\,\bigr\}.
  \end{split} 
  \end{equation*}
\end{definition}

Note that $\seg{i}{\xi_2} = \seg{i+1}{\xi_1}$, $i \in\{ 1,\ldots,\nseg-1\}$,
and $\seg{1}{\xi_1} \in \RR$ is an arbitrary origin for the composed interval $\Omega$.
For all practical purposes, hereafter we will assume $\seg{i}{\contin} \leq \min (\seg{i}{\pdeg},\seg{i+1}{\pdeg})$, where $\seg{i}{\pdeg}$ stands for the degree of the spline space defined on the $i$-th segment.
In the special case of $\seg{i}{\contin} = \seg{i}{\pdeg}=\seg{i+1}{\pdeg}$ we get a $C^\infty$ smooth segment join.

It is easy to see that increasing the continuity orders $\mbf{\contin}$ at the segment joins 
leads to smoother subspaces. 
\begin{proposition} \label{pro:nestedness}
  If $\mbf{\contin_1}$ and $\mbf{\contin_2}$ are such that $\seg{i}{\contin_1} \leq \seg{i}{\contin_2}$ for all $i$, 
  then $\SPL{S}{\contin_1}{U} \supseteq \SPL{S}{\contin_2}{U}$.
\end{proposition}

We can also reformulate the standard (uniform degree) spline spaces $\BSP{}{U}$ in terms of the above definition, starting from local polynomial spaces (see Remark~\ref{rmk:bezier}).

\begin{proposition} \label{pro:bezier}
Given $\pdeg\in\NN$, $\mbf{\intIx} := (\seg{1}{\intIx},\seg{2}{\intIx},\ldots,\seg{\nseg+1}{\intIx})\in \RR^{\nseg+1}$ with $\seg{i}{x}<\seg{i+1}{x}$, and $\mbf{\contin} := (\seg{1}{\contin},\seg{2}{\contin},\ldots,\seg{\nseg-1}{\contin})\in \NN_{-1}^{\nseg-1}$, we set 
\begin{equation*}
  \mbf{U} := (\seg{1}{U},\dots,\seg{\nseg}{U}), \quad
  \seg{i}{U} := [\;\underbrace{\seg{i}{\intIx},\ldots,\seg{i}{\intIx}}_{\pdeg+1 \text{ times}}\;,
  \;\underbrace{\seg{i+1}{\intIx},\ldots,\seg{i+1}{\intIx}}_{\pdeg+1 \text{ times}}\;], \quad 
  i=1,\ldots,\nseg,
\end{equation*}
and
\begin{equation*}
  U := [\;\underbrace{\seg{1}{\intIx},\ldots,\seg{1}{\intIx}}_{\pdeg+1 \text{ times}}\;, 
  \;\underbrace{\seg{2}{\intIx},\ldots,\seg{2}{\intIx}}_{\pdeg-\seg{1}{\contin} \text{ times}}\;,
  \;\ldots\;,
  \;\underbrace{\seg{\nseg}{\intIx},\ldots,\seg{\nseg}{\intIx}}_{\pdeg-\seg{\nseg-1}{\contin} \text{ times}}\;,
  \;\underbrace{\seg{\nseg+1}{\intIx},\ldots,\seg{\nseg+1}{\intIx}}_{\pdeg+1 \text{ times}}\;].
\end{equation*}
Then,
\begin{equation*}
  \SPL{S}{\contin}{U}=\BSP{}{U}.
\end{equation*}
\end{proposition}

\begin{remark} \label{rmk:polynomial-vs-spline}
Definition~\ref{def:multiDegSpline} takes as input a set of open knot vectors $\seg{i}{U}$, $i=1,\ldots,\nseg$ (specifying the uniform degree spline spaces $\BSP{i}{U}$ as explained in Remark~\ref{rmk:structureU}) and a set of continuity orders $\seg{i}{\contin}$, $i=1,\ldots,\nseg-1$ (specifying how smoothly to glue together the given spline spaces). 
This setup differs from the multi-degree spline definition in \cite{beccari2017} where polynomial spaces are glued together with a certain smoothness. Even though both definitions are conceptually equivalent, Definition~\ref{def:multiDegSpline} is 
computationally more interesting
because it enables us to merge consecutive polynomial spaces of the same degree into a single (uniform degree) spline space; see Proposition~\ref{pro:bezier}. This reduces the number of segment joins where a special gluing technique needs to be applied and leads to more efficient constructions of multi-degree spline spaces and basis functions.
\end{remark}

\begin{remark}
A more general multi-degree spline framework is considered in \cite{toshniwal2017polar}, allowing for local rational spline (NURBS) spaces to be joined smoothly, and also incorporating periodic smoothness constraints. We do not consider this full generality here, in order to focus on the essentials of the basis construction. A discussion on how to extend the implementation to the periodic case is provided in Remark~\ref{rmk:periodic}.
\end{remark}

In the following section we describe a procedure to construct a basis for the multi-degree spline space $\SPL{S}{\contin}{U}$. We do this by formulating the continuity constraints that splines locally defined over adjacent segments would need to satisfy. The (left) null-space of these constraints yields the coefficients of an extraction operator that help us to create the multi-degree basis functions as linear combinations of the local B-splines on each segment.

\subsection{Multi-degree spline extraction operator}

Consider an $\nseg$-segment knot vector configuration $\mbf{U}$.
On the $i$-th knot vector $\seg{i}{U}$, we have $\seg{i}{\ndof}$ B-splines $\seg{i}{\bs_{j,\seg{i}{\pdeg}}}$ of degree $\seg{i}{\pdeg}$ that span the spline space $\BSP{i}{U}$. 
In the first step, we map these basis functions from $\seg{i}{\intI}$ to $\seg{i}{\Omega}$ using $\seg{i}{\phi}$ in Definition \ref{def:multiDegSpline}, and extend them on the entire interval $\Omega$ by defining them to be zero outside $\seg{i}{\Omega}$.
More precisely, setting 
\begin{equation}\label{eq:dim-cumsum}
  \nsum_{0}:=0, \quad 
  \nsum_i:=\sum_{\ell=1}^{i}\seg{\ell}{\ndof}=\nsum_{i-1}+\seg{i}{\ndof},\quad i=1,\ldots,\nseg, 
\end{equation}
we define for $i\in\{1,\ldots,n_s\}$, $j\in\{1,\ldots,\seg{i}{\ndof}\}$,
\begin{equation}\label{eq:basis-start}
  \bsext_{\nsum_{i-1} + j}(\xi) := 
  \begin{cases}
  \; \seg{i}{\bs_{j,\seg{i}{\pdeg}}}(\intIx),& 
  \text{if } [\seg{i}{\xi_1},\;\seg{i}{\xi_2}) \ni \xi = \seg{i}{\phi}(\intIx),\\[0.2cm]
  \; \seg{i}{\bs_{j,\seg{i}{\pdeg}}}(\seg{i}{\intIx_2}),& 
  \text{if } i=\nseg \text{ and } \xi=\seg{\nseg}{\xi_2},\\[0.2cm]
  \; 0, & \text{otherwise}.
  \end{cases}
\end{equation}
For the sake of simplicity, we dropped the reference to the (local) degree in the notation.
From the properties of B-splines, it is clear that the functions $\bsext_{1},\ldots,\bsext_{\nsum_{\nseg}}$ are linearly independent,
form a non-negative partition of unity, and span the space $\SPL{S}{-1}{U}$, where $\mbf{-1}=(-1,\ldots,-1)$.
We arrange these basis functions in a column vector $\mbf{\bsext}$ of length $\nsum_{\nseg}$.

Now, we are looking for a set of basis functions $\{\mdbs_1,\ldots,\mdbs_\ndof\}$ that span the smoother space $\SPL{S}{\contin}{U}$.
We arrange these basis functions in a column vector $\mbf{\mdbs}$ of length $\ndof$.
From Proposition~\ref{pro:nestedness} we know that $\SPL{S}{\contin}{U}\subseteq\SPL{S}{-1}{U}$.
Hence, we aim to construct a matrix $\mbf{H}$ of size $n \times \nsum_{\nseg}$ such that 
\begin{equation}\label{eq:multiDegBasisExt}
  \mbf{\mdbs} = \mbf{H}\,\mbf{\bsext}.  
\end{equation}
To this end, we build continuity constraints at all segment joins corresponding to $\mbf{\contin}$ and construct $\mbf{H}$ as their (left) null-space. 

\subsubsection{Continuity constraint matrix}\label{sec:contConst}
Consider the $i$-th segment join for some $i\in\{1,\ldots,\nseg-1\}$.
Let $\seg{i,1}{\mbf{K}}$ be a matrix of size $\seg{i}{\ndof}\times(\seg{i}{\contin}+1)$,
whose $\deriv$-th column is given by
\begin{equation} \label{eq:K_L}
    \begin{bmatrix}
    0 & \cdots & 0 & D_-^{\deriv-1}\bsext_{\nsum_i-\deriv+1}(\seg{i}{\xi_2}) & \cdots & D_-^{\deriv-1}\bsext_{\nsum_i}(\seg{i}{\xi_2})
    \end{bmatrix}^T,
\end{equation}
and let $\seg{i,2}{\mbf{K}}$ be a matrix of size $\seg{i+1}{\ndof}\times(\seg{i}{\contin}+1)$, whose $\deriv$-th column is given by
\begin{equation} \label{eq:K_R}
    \begin{bmatrix}
    -D_+^{\deriv-1}\bsext_{\nsum_i+1}(\seg{i+1}{\xi_1}) & \cdots & -D_+^{\deriv-1}\bsext_{\nsum_i+\deriv}(\seg{i+1}{\xi_1}) & 0 & \cdots & 0
    \end{bmatrix}^T.
\end{equation}
Note that the derivatives of the basis functions in the above matrices can be easily computed by evaluating the derivatives of the corresponding local B-splines at the end points of their basic interval. This explains the triangular structure of both matrices.

Using these matrices, we can build the matrix $\seg{i}{\mbf{K}}$ of size $\nsum_{\nseg}\times(\seg{i}{\contin}+1)$ which contains all constraints required to enforce $C^{\seg{i}{\contin}}$ at $\seg{i}{\xi_2}=\seg{i+1}{\xi_1}$. This matrix is defined row-wise in the following manner:
\begin{itemize}
  \item the $(\nsum_{i-1}+\ell)$-th row of $\seg{i}{\mbf{K}}$ is equal to the $\ell$-th row of $\seg{i,1}{\mbf{K}}$;
  \item the $(\nsum_i+\ell)$-th row of $\seg{i}{\mbf{K}}$ is equal to the $\ell$-th row of $\seg{i,2}{\mbf{K}}$;
  \item all other rows of $\seg{i}{\mbf{K}}$ are identically zero.
\end{itemize}
It can be easily verified that for a row vector of coefficients $\mbf{f}$ such that $\mbf{f}\,\seg{i}{\mbf{K}} = \mbf{0}$, the spline defined by $\mbf{f}\,\mbf{\bsext}$ is going to be $C^{\seg{i}{\contin}}$ across $\seg{i}{\xi_2}$. Therefore, once all the matrices $\seg{i}{\mbf{K}}$ have been assembled, the only remaining step is the construction of $\mbf{H}$ such that it spans their left null-spaces. The matrix $\mbf{H}$ is called \emph{multi-degree spline extraction operator}, and we employ the algorithm in Figure \ref{alg:extraction} for its construction.

\begin{figure}[t!]
\begin{algorithm}[H]
  \small
  \Fn{extraction}{
  \In{continuity constraint matrices $\seg{i}{\mbf{K}}$ (size: $\nsum_\nseg \times (\seg{i}{\contin}+1)$) for $i=1,\ldots,\nseg-1$}
  \Out{extraction matrix $\mbf{H}$}
  \bigskip
  $\mbf{H} \gets \mbox{identity matrix}$ (size: $\nsum_\nseg \times \nsum_\nseg$)\;
  \For{$i \gets 1$ \KwTo $\nseg-1$}{
    $\mbf{L} \gets \mbf{H}\, \seg{i}{\mbf{K}}$\;
    \For{$j \gets 0$ \KwTo $\seg{i}{\contin}$}{
      $\bar{\mbf{H}} \gets \mbox{sparse null-space of $(j+1)$-th column of } \mbf{L}$\;
      $\mbf{H} \gets \bar{\mbf{H}}\, \mbf{H}$\;
      $\mbf{L} \gets \bar{\mbf{H}}\, \mbf{L}$\;
    }
  }
  }
\end{algorithm}
\caption{Computation of extraction matrix $\mbf{H}$ in pseudo-code}\label{alg:extraction}
\end{figure}

Essentially, the algorithm works as follows. The outer for-loop runs over all the segment joins, and the inner for-loop over all the 
continuity constraints at each join. We address a single continuity constraint at a time, so we increase the smoothness of the basis functions (obtained by \eqref{eq:multiDegBasisExt}) gradually.
The matrix $\mbf{L}$ keeps track of the remaining continuity constraints for the basis functions built so far.
In each step, we first construct the left null-space $\bar{\mbf{H}}$ of the next continuity constraint in the matrix $\mbf{L}$ (detailed in Section \ref{sec:nullSpace}). 
Then, we update the global null-space matrix $\mbf{H}$ by left multiplying it with $\bar{\mbf{H}}$, thus incorporating the current continuity constraint.
Finally, we update $\mbf{L}$ by again left multiplying it with $\bar{\mbf{H}}$, thus obtaining a reduced system of constraints that need to be satisfied. 
The process is repeated until all constraints at all segment joins have been satisfied.

\begin{remark}
There are two main differences between the algorithm in Figure \ref{alg:extraction} and the original algorithm presented in \cite[Algorithm~1]{toshniwal2017polar}. 
\begin{itemize}
  \item The outer and inner for-loop are swapped. This swap has three beneficial consequences. 
  First, it may reduce the number of steps in the present inner for-loop. There are now $\seg{i}{\contin}+1$ steps instead of $\max_\ell\seg{\ell}{\contin}+1$; hence, it is more efficient when requiring different orders of continuity at the different segment joins.
  Second, we do not need to store all the matrices $\seg{i}{\mbf{K}}$, $i=1,\ldots,\nseg-1$ in memory, but we can work with a single matrix $\mbf{L}$. 
  Third, we can avoid an additional for-loop with matrix multiplications updating the different matrices $\seg{i}{\mbf{K}}$ in each step; it suffices to do a single matrix multiplication at the initialization of $\mbf{L}$.
  \item We compute the left null-space instead of the (right) null-space. This eliminates a matrix transpose in each step of the inner for-loop.
\end{itemize}

\end{remark}

\subsubsection{Sparse null-space construction}\label{sec:nullSpace}

We now focus on the left null-space computation of a column of the continuity constraint matrix $\mbf{L}$.
Any basis of the null-space would lead to a valid basis of the space $\SPL{S}{\contin}{U}$ using the previously described procedure. However, we are not just interested in any basis, but are looking for a basis with similar properties to the standard B-spline basis. In particular, we embrace properties like local support and non-negative partition of unity. This brings us to the concept of IGA-suitable extraction operator \cite{toshniwal2017polar}.

\begin{definition}[IGA-suitable extraction] \label{def:dasExtraction}
  An extraction operator $\mbf{H}$, acting on locally supported basis functions that form a non-negative partition of unity, is called IGA-suitable if
  \begin{itemize}
    \item $\mbf{H}$ is a full-rank matrix,
    \item each column of $\mbf{H}$ sums to 1,
    \item each entry in $\mbf{H}$ is non-negative, and
    \item $\mbf{H}$ has a band-sparsity pattern.
  \end{itemize}  
\end{definition}
It is easy to verify that the action of an IGA-suitable extraction operator on a local basis forming a non-negative partition of unity gives rise to another local basis that also forms a non-negative partition of unity. Keeping in mind the properties of the starting basis \eqref{eq:basis-start}, we aim for the matrix $\mbf{H}$ in \eqref{eq:multiDegBasisExt} to be IGA-suitable. This is obtained when all intermediate matrices $\bar{\mbf{H}}$ in the procedure possess the same feature. We employ the algorithm in Figure \ref{alg:nullSpace} for its construction.

\begin{figure}[t!]
\begin{algorithm}[H]
  \small
  \Fn{nullspace}{
  \In{vector $\mbf{l}$ (length $q$)}
  \Out{left null-space matrix $\bar{\mbf{H}}$}
  \bigskip
  $\bar{\mbf{H}} \gets \mbox{zero matrix}$ (size: $(q-1) \times q$)\;
  $i_1 \gets$ index of first non-zero element of $\mbf{l}$\;
  $i_2 \gets$ index of last non-zero element of $\mbf{l}$\;
  \For{$j \gets 1$ \KwTo $i_1$}{
    $\bar{\mbf{H}}(j,j) \gets 1$\;
  }
  \For{$j \gets i_1+1$ \KwTo $i_2-1$}{
    $\bar{\mbf{H}}(j-1,j) \gets -\dfrac{\mbf{l}(j-1)}{\mbf{l}(j)}\,\bar{\mbf{H}}(j-1,j-1)$\;
    $\bar{\mbf{H}}(j,j) \gets 1 - \bar{\mbf{H}}(j-1,j)$\;
  }
  \For{$j \gets i_2$ \KwTo $q$}{
    $\bar{\mbf{H}}(j-1,j) \gets 1$\;
  }
  }
\end{algorithm}
\caption{Computation of sparse null-space $\bar{\mbf{H}}$ of vector $\mbf{l}$ (column of $\mbf{L}$) in pseudo-code}\label{alg:nullSpace}
\end{figure}

Essentially, we strive to build the sparsest possible left null-space of column vector $\mbf{l}$, containing the next continuity constraint in the matrix $\mbf{L}$. Because of the triangular structure of the matrices in \eqref{eq:K_L}--\eqref{eq:K_R} and the B-spline properties, such vector $\mbf{l}$ has a special sparsity structure: all its entries are zero, except for a single, consecutive block of non-zero entries whose sum equals zero.
Assuming $\mbf{l}=[l_1,\ldots,l_q]^T$, we build the matrix $\bar{\mbf{H}}$ as follows. For $i=1,\ldots,q$,
\begin{itemize}
  \item if $l_i=0$, then we add the row
\begin{equation*}
  [0,\;\dots,\;0,\;\underbrace{\;\;\;\; 1 \;\;\;\;}_{i\text{-th entry}}\;,\;0,\;\dots,\;0];
\end{equation*}
  \item if $l_{i}\neq0$ and $l_{i+1}\neq0$, then we add the row
\begin{equation*}
  [0,\;\dots,\;0,\;\underbrace{\; 1-h\vphantom{\frac{l_i}{l_i}}\;}_{i\text{-th entry}}\;,\;\underbrace{\; -\frac{l_{i}}{l_{i+1}}(1-h) \;}_{(i+1)\text{-th entry}}\;,\;0,\;\dots,\;0],
\end{equation*}
where $h$ is the value at the $i$-th entry of the previous row when $i\geq2$, or zero otherwise.
\end{itemize}
Assuming that the non-zero entries of $\mbf{l}$ are indexed from $i_1$ to $i_2$, we deduce from their summation-to-zero that
\begin{equation*}
-\frac{l_{i_2-1}}{l_{i_2}}\left(1+\frac{l_{i_2-2}}{l_{i_2-1}}\left(1+\cdots \frac{l_{i_1+2}}{l_{i_1+3}}\left(1+\frac{l_{i_1+1}}{l_{i_1+2}}\left(1+\frac{l_{i_1}}{l_{i_1+1}}\right)\right)\cdots\right)\right)=1.
\end{equation*}
From this identity and the above construction, we directly conclude that each row of $\bar{\mbf{H}}$ is in the left null-space of $\mbf{l}$ and that each column of $\bar{\mbf{H}}$ sums to 1.
Furthermore, we see that $\bar{\mbf{H}}$ is an upper bidiagonal matrix. This property can be exploited in an efficient implementation (see Remark~\ref{rmk:matlab}). 

\begin{figure}[t!]
\begin{algorithm}[H]
  \small
  \Fn{nullspace (\Matlab implementation)}{
  \In{\Matlab vector \texttt{ll}}
  \Out{\Matlab matrix {\texttt{Hbar}} representing the left null-space}
  \bigskip
  \texttt{q = length(ll)}\; 
  \texttt{i1 = find(ll, 1, `first')}\; 
  \texttt{i2 = find(ll, 1, `last')}\;
  \texttt{dd = zeros(q-1, 2)}\;
  \texttt{dd(1:i1, 1) = 1}\;
  \ForM{\texttt{j = i1:i2-2}}{
    \texttt{dd(j, 2) = -ll(j) / ll(j+1) * dd(j, 1)}\;
    \texttt{dd(j+1, 1) = 1 - dd(j, 2)}\;
  }
  \texttt{dd(i2-1:q-1, 2) = 1}\;
  \texttt{Hbar = spdiags(dd, [0 1], q-1, q)}\;
  }
\end{algorithm}
\caption{Efficient computation of sparse null-space $\texttt{Hbar}$ of vector $\texttt{ll}$  (column of $\mbf{L}$) in \Matlab}\label{alg:nullSpace-matlab}
\end{figure}

\begin{remark}\label{rmk:conjecture}
  It has been conjectured in \cite{toshniwal2017polar} that the matrix $\bar{\mbf{H}}$ contains only non-negative entries,
  and forms an IGA-suitable extraction operator. This was formally proven in \cite{toshniwal2018mdb}.
\end{remark}

\begin{remark}\label{rmk:matlab}
  The described extraction mechanism can be efficiently encoded in \Matlabs. 
  Since both the extraction matrix $\mbf{H}$ and each intermediate null-space matrix $\bar{\mbf{H}}$ are sparse, it is natural to store them as \Matlab sparse matrices. In particular, since $\bar{\mbf{H}}$ is an upper bidiagonal matrix, it can be constructed in a fast and memory-efficient manner via the \Matlab command \texttt{spdiags}. The first and last for-loop in Figure~\ref{alg:nullSpace} can also be avoided by directly using a sub-matrix assignment in \Matlabs. For the sake of completeness, an efficient \Matlab implementation is depicted in Figure~\ref{alg:nullSpace-matlab}. It is clear that both its computational and memory complexity are linear with respect to the length of the input vector. This implies that the null-space construction has an optimal complexity.
\end{remark}
%

\begin{remark}
  The original null-space construction described in \cite[Section~2.3.2]{toshniwal2017polar} has a rather cumbersome procedure to select the indices of non-zero entries in $\mbf{l}$. This is mainly due to the fact that periodically smooth configurations are directly incorporated into the computation of the extraction matrix. 
  By explicitly avoiding the periodic case, the index selection can be implemented in an elegant way as illustrated in Figure~\ref{alg:nullSpace-matlab}, using the variables \texttt{i1} and \texttt{i2}.
\end{remark}

\begin{remark}\label{rmk:periodic}
  Even though the algorithms presented in Figures~\ref{alg:extraction}--\ref{alg:nullSpace} are not explicitly designed for building multi-degree spline functions with a periodic nature 
  (meaning that the derivatives at both end points of the domain are equal up to a certain order),
  they can still be useful in that context with a minimal extra effort and under certain mild restrictions. There are several possible strategies to this aim. 
  A first strategy is to consider a ``periodic extension'' of the spline segments (i.e., copying few spline segments from one end to the other end), to apply the algorithms in Figures~\ref{alg:extraction}--\ref{alg:nullSpace}, and then to reassemble the extended non-periodic extraction matrix into the periodic one. 
  An alternative strategy is to circularly shift the rows of the extraction matrix built so far, such that the periodic continuity constraints behave like continuity constraints at an interior segment join, and then the null-space computation of a column of the new matrix $\mbf{L}$ can be simply done by the algorithm in Figure~\ref{alg:nullSpace}. 
  This can be efficiently encoded in \Matlab using the command \texttt{circshift} and has been implemented 
  in the accompanying \Matlab toolbox.
\end{remark}

\subsection{Properties of multi-degree splines}\label{sss:multiDegSplineEx}

The spline functions $\mdbs_1,\ldots,\mdbs_\ndof$ constructed by means of the multi-degree spline extraction operator $\mbf{H}$ in \eqref{eq:multiDegBasisExt} are called \emph{MDB-splines}. Several interesting properties of MDB-splines were proven in \cite{toshniwal2017polar,toshniwal2018mdb}; we just give a selection here.

\begin{proposition}\label{pro:properties-mdbspline}
The MDB-splines $\mdbs_j$, $j=1,\ldots,\ndof$ enjoy the following properties:
\begin{itemize}
  \item local support,
  \item $C^{\seg{i}{\contin}}$-smoothness at the $i$-th segment join,
  \item non-negative partition of unity, and
  \item interpolation at the end points.
\end{itemize}
\end{proposition}
%
Furthermore, MDB-splines are linearly independent, so the multi-degree spline space $\SPL{S}{\contin}{U}$ has dimension
\begin{equation} \label{eq:dimension-md}
  n = \seg{1}{\ndof} + \sum_{i=2}^{\nseg} (\seg{i}{\ndof} - \seg{i-1}{\contin} - 1).
\end{equation}

Overall, MDB-splines behave very similarly to B-splines (cf. Proposition~\ref{pro:properties-bspline}). Actually, if all degrees are chosen equal, i.e., $\seg{1}{\pdeg}=\cdots=\seg{\nseg}{\pdeg}=\pdeg$, then the resulting MDB-splines are exactly the B-splines of degree $\pdeg$ forming the same space $\SPL{S}{\contin}{U}$.

We can represent a multi-degree spline $s\in\SPL{S}{\contin}{U}$ in terms of MDB-splines using a vector of coefficients $\mbf{s} := [s_1,\dots,s_\ndof]$. Then,
\begin{equation} \label{eq:multiDegFun}
s(\xi) = \sum_{j=1}^\ndof s_j\mdbs_j(\xi) = \sum_{j=1}^\ndof s_j \sum_{i=1}^{\nsum_{\nseg}}H_{ji}\bsext_i(\xi)
= \sum_{i=1}^{\nsum_{\nseg}} f_i\bsext_i(\xi), \quad \xi\in\Omega,
\end{equation}
where $H_{ji}$ is the $(j,i)$-th entry of $\mbf{H}$, and
\begin{equation*} 
f_i := \sum_{j=1}^\ndof s_j H_{ji}.
\end{equation*}
The latter coefficients can be compactly written in vector notation as
\begin{equation*}
\mbf{f}:=\mbf{s}\,\mbf{H}.
\end{equation*}

\begin{remark} \label{rmk:from-mdb-to-b}
The multi-degree spline extraction operator allows us to move directly from the global MDB-spline representation to the local B-spline representation, and then apply any of the existing efficient B-spline algorithms. This delegation avoids that new algorithms have to be developed for each kind of MDB-spline manipulation (evaluation, differentiation, etc.), which does not only simplify implementation but also profit optimally of the powerful B-spline technology. 
As a consequence, when precomputing and storing the corresponding (sparse) extraction matrix, the efficiency and complexity of such algorithms to manipulate MDB-splines is entirely in the hands of the considered B-spline algorithms (and their implementations).
\end{remark}

\begin{remark} \label{rmk:bezier-extraction}
If the multi-degree spline extraction operator $\mbf{H}$ is computed for a spline space as in Proposition~\ref{pro:bezier}, then it is nothing other than the so-called \emph{B\'ezier extraction operator}. Such an operator provides the conversion from the global (MD)B-spline representation to the local B\'ezier representation (see Remark~\ref{rmk:bezier}).
B\'ezier extraction is a very common procedure in IGA to fit spline representations in the standard finite element pipeline (see, for example, \cite{borden2011isogeometric}). Afterwards, the B\'ezier representation could then easily be converted into an interpolatory Lagrange representation, if wanted (see \cite{schillinger2016}).
\end{remark}

Refinement can also be performed on any multi-degree spline in $\SPL{S}{\contin}{U}$ by exploiting the multi-degree spline extraction operator. We refer the reader to \cite[Section~2.4.3]{toshniwal2017polar} for details about (local) mesh refinement and (local) degree raising on one or more of its segments. 
In summary, a spline $s\in\SPL{S}{\contin}{U}$ represented as in \eqref{eq:multiDegFun} can be converted into the (locally) refined space $\SPL{S}{\tilde{\contin}}{\tilde{U}}\supset\SPL{S}{\contin}{U}$ as follows, under the assumption that $\mbf{U}$ and $\mbf{\tilde{U}}$ have the same number of B-spline segments. Given
\begin{itemize}
  \item the old and new extraction operators: $\mbf{H}$ and $\mbf{\tilde{H}}$,
  \item for each segment $i=1,\ldots,\nseg$ the local B-spline refinement matrix $\seg{i}{\mbf{R}}$ from 
  $\BSP{i}{U}$ to $\BSP{i}{\tilde{U}}$,
  converting any vector of B-spline coefficients $\seg{i}{\mbf{f}}$ into the new vector of B-spline coefficients
  $\seg{i}{\mbf{\tilde{f}}}=\seg{i}{\mbf{f}}\,\seg{i}{\mbf{R}}$,
\end{itemize}
the new vector of MDB-spline coefficients $\mbf{\tilde{s}}$ can be computed as
\begin{equation} \label{eq:conversion-mdb}
\mbf{\tilde{s}} = \mbf{s} \,\mbf{H}\,\diag(\seg{1}{\mbf{R}},\ldots,\seg{\nseg}{\mbf{R}})\,
\mbf{\tilde{H}}^T (\mbf{\tilde{H}}\,\mbf{\tilde{H}}^T)^{-1}.
\end{equation}

\begin{remark}
The refined coefficients in \eqref{eq:conversion-mdb} are basically obtained via a least-squares procedure. 
In particular, the multiplication on the right 
with
$\mbf{\tilde{H}}^T (\mbf{\tilde{H}}\,\mbf{\tilde{H}}^T)^{-1}$ can be efficiently encoded in \Matlab using the (forward) slash operation with $\mbf{\tilde{H}}$.
All the matrices involved are (highly) sparse, so the whole computation is not so costly. This approach has been adopted in the accompanying \Matlab toolbox for handling conversions between different multi-degree spline spaces.
\end{remark}

\begin{remark}
The conversion procedure in \eqref{eq:conversion-mdb} is quite abstract and hence very generally applicable. 
It shows us how any kind of refinement in the MDB-spline setting can be delegated to the corresponding refinement in the local B-spline setting (the matrices $\seg{i}{\mbf{R}}$). This does not only simplify implementation but also profit optimally of the known B-spline technology. Hence, the efficiency and complexity of such algorithms for (locally) refining MDB-splines strongly depends on the considered B-spline algorithms (see also Remark~\ref{rmk:from-mdb-to-b}).
\end{remark}

\begin{remark}
In the recent literature there exist alternative ways to perform mesh refinement or degree raising on multi-degree splines. For example, specific refinement procedures based on so-called transition functions are described in \cite[Section~5]{beccari2017}. 
\end{remark}

\section{Examples} \label{sec:matlab}

We conclude with three basic examples to illustrate the construction and use of MDB-splines through the previously described extraction procedure implemented in a small \Matlab toolbox.  
Full details of the facilities available from the \Matlab toolbox may be found in the user manual 
that accompanies the software which is available through the CALGO library.

In the first and second example we construct and visualize some MDB-splines (see the test code \texttt{EX\_basis.m}) and periodic MDB-splines (see the test code \texttt{EX\_basis\_periodic.m}). In the third example we compare the MDB-spline representation with the standard B-spline representation 
(see the test code \texttt{EX\_conversion.m}).

\begin{figure}[t!]
\subfigure[$\contin=0$]{\includegraphics[width=4.4cm]{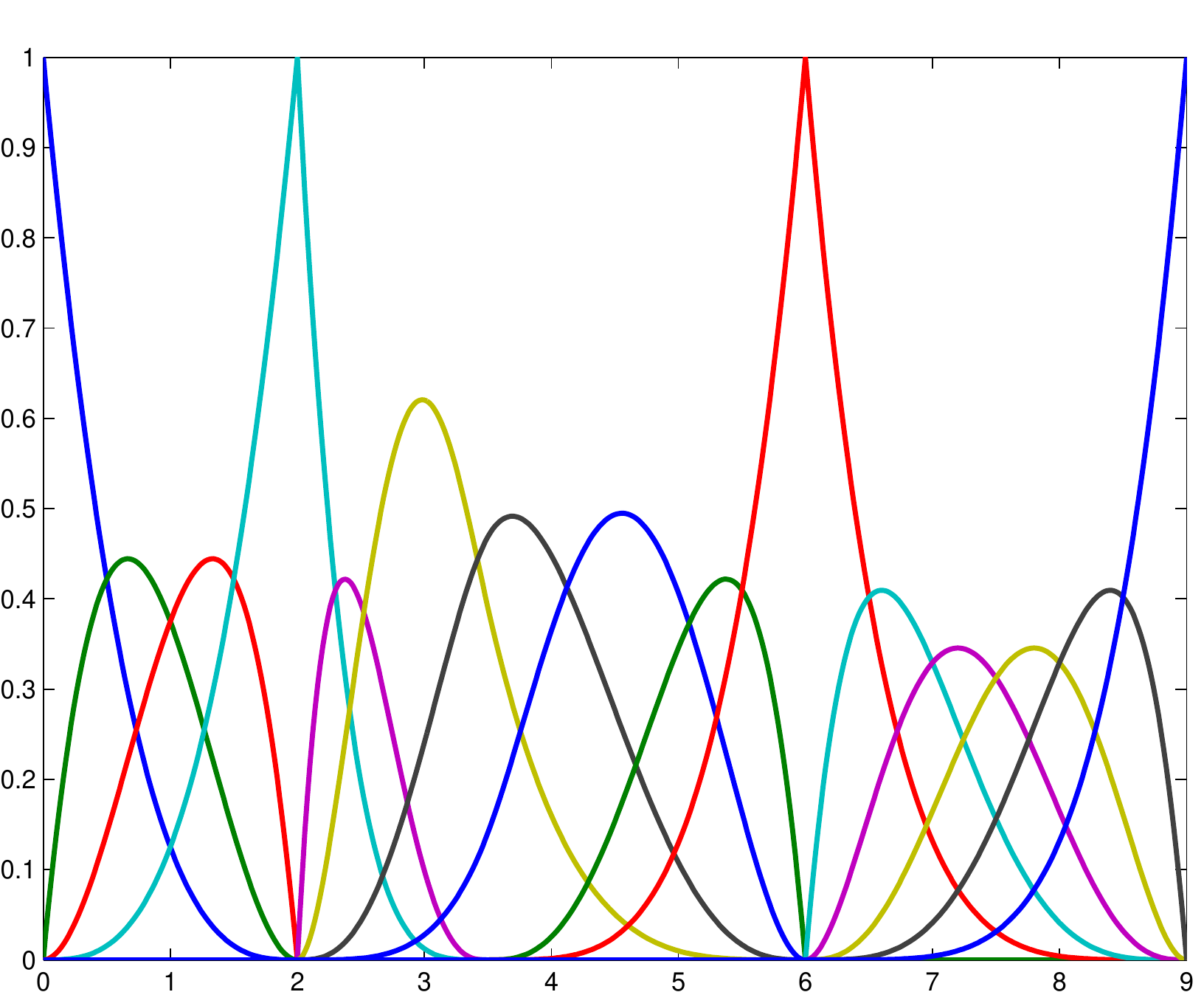}} \hspace*{0.1cm}
\subfigure[$\contin=1$]{\includegraphics[width=4.4cm]{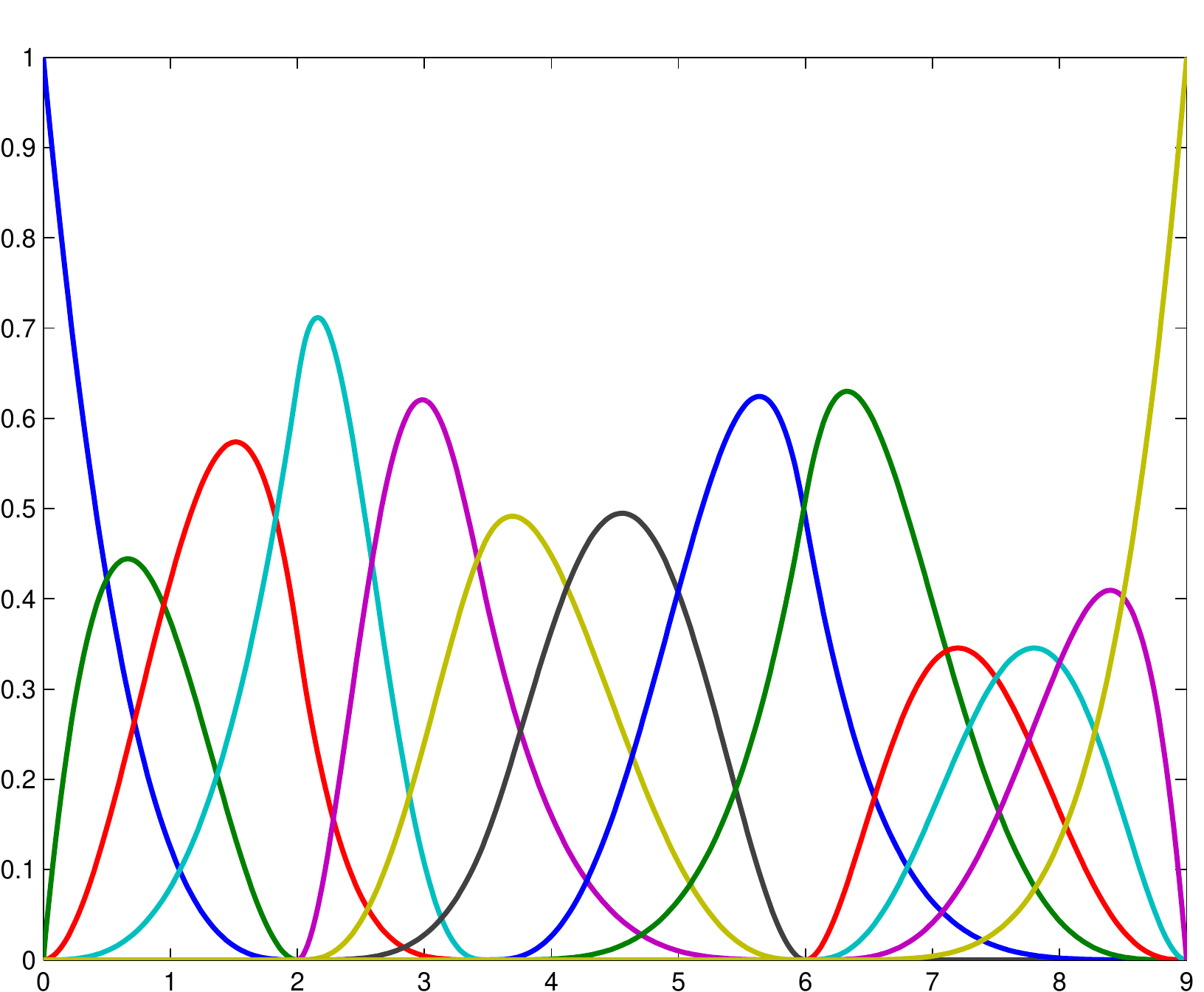}} \hspace*{0.1cm}
\subfigure[$\contin=2$]{\includegraphics[width=4.4cm]{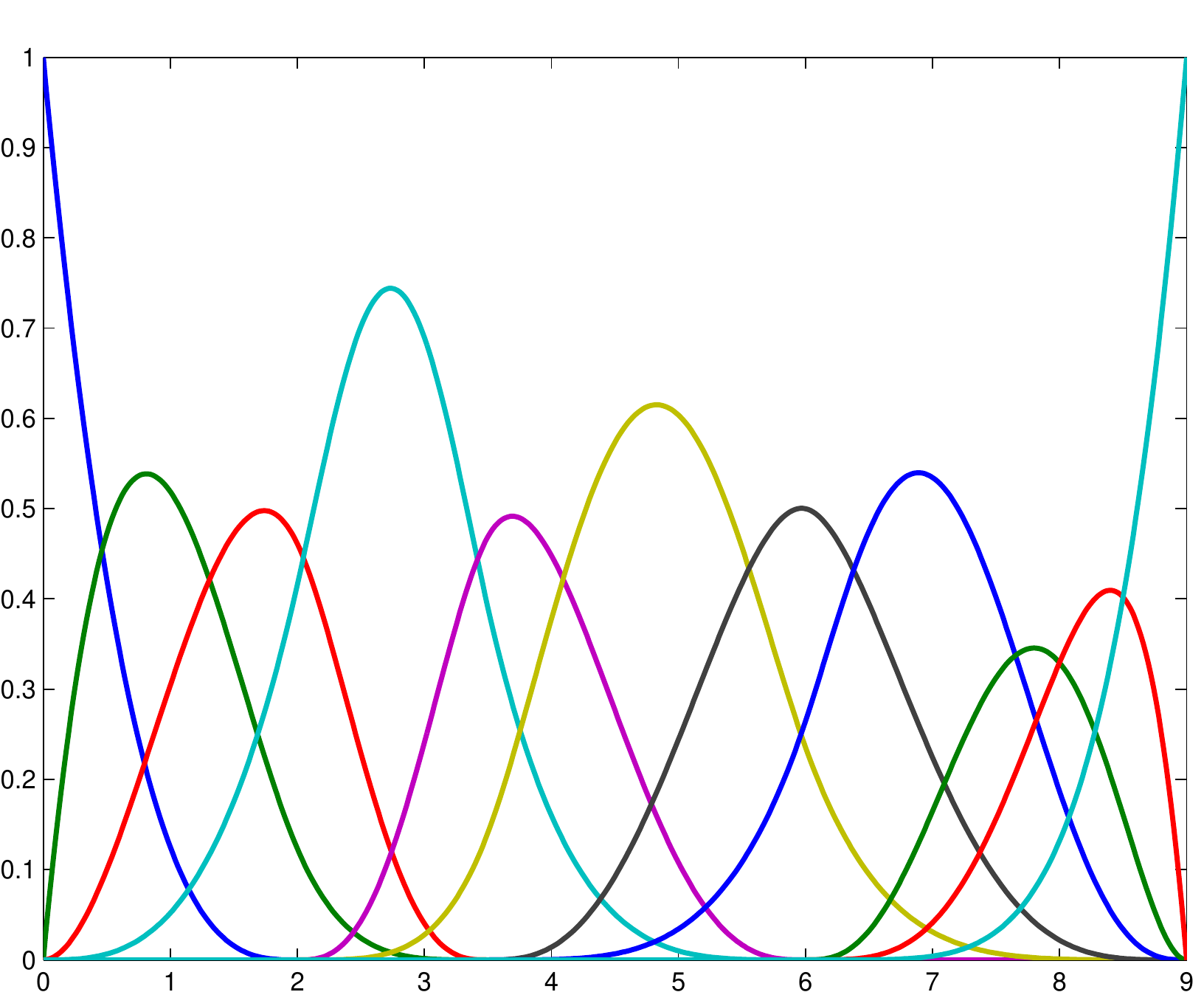}} \vspace*{-0.2cm}
\caption{MDB-spline basis of degrees $(3,4,5)$ defined over the knot vector configuration $\mbf{U}=(\seg{1}{U},\seg{2}{U},\seg{3}{U})$ specified in \eqref{eq:ex1-knots} and continuity orders $\mbf{\contin}=(\contin,\contin)$ for $\contin\in\{0,1,2\}$ }
\label{fig:ex-basis}
\bigskip
\subfigure[$\contin=0$]{\includegraphics[width=4.4cm]{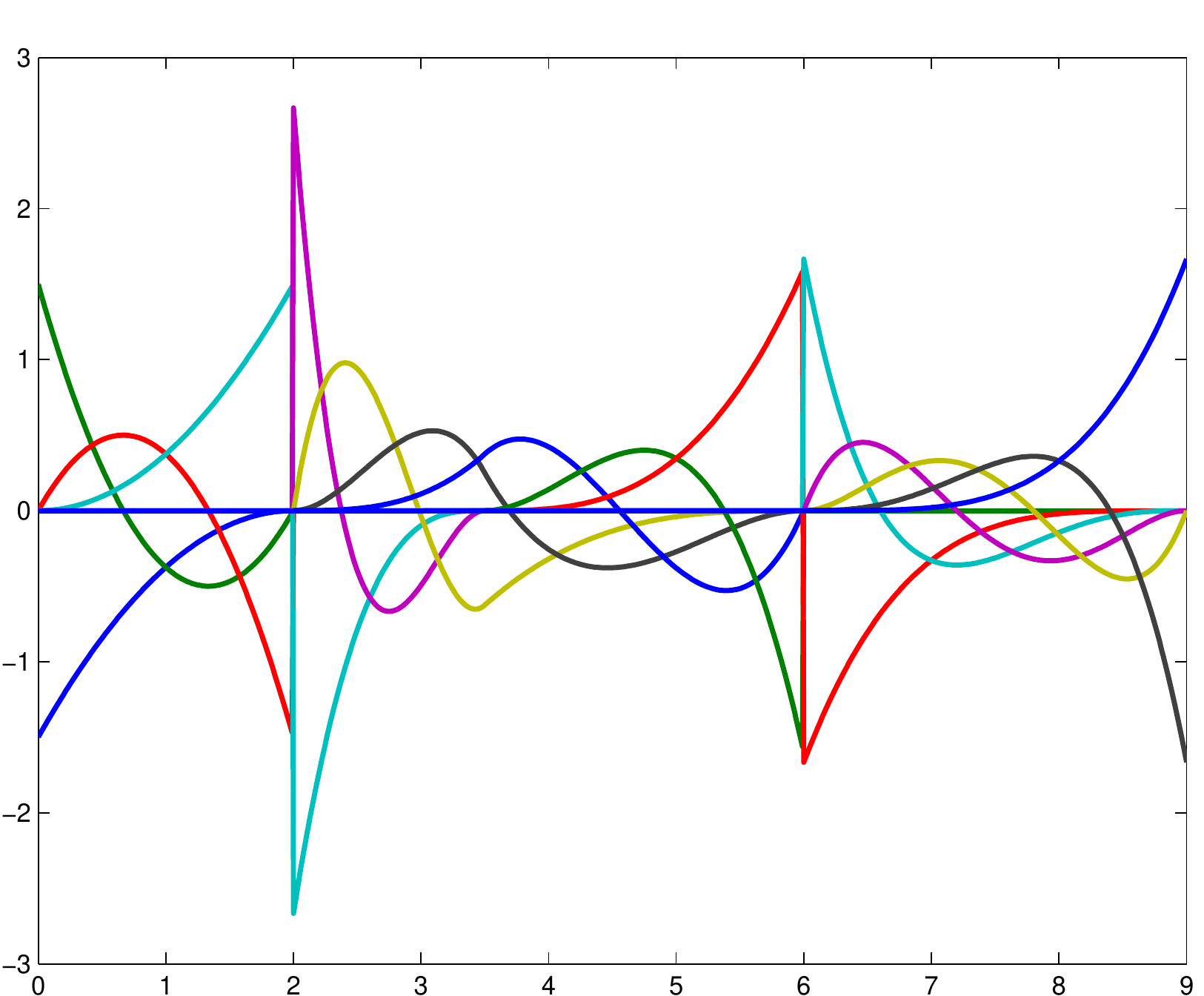}} \hspace*{0.1cm}
\subfigure[$\contin=1$]{\includegraphics[width=4.4cm]{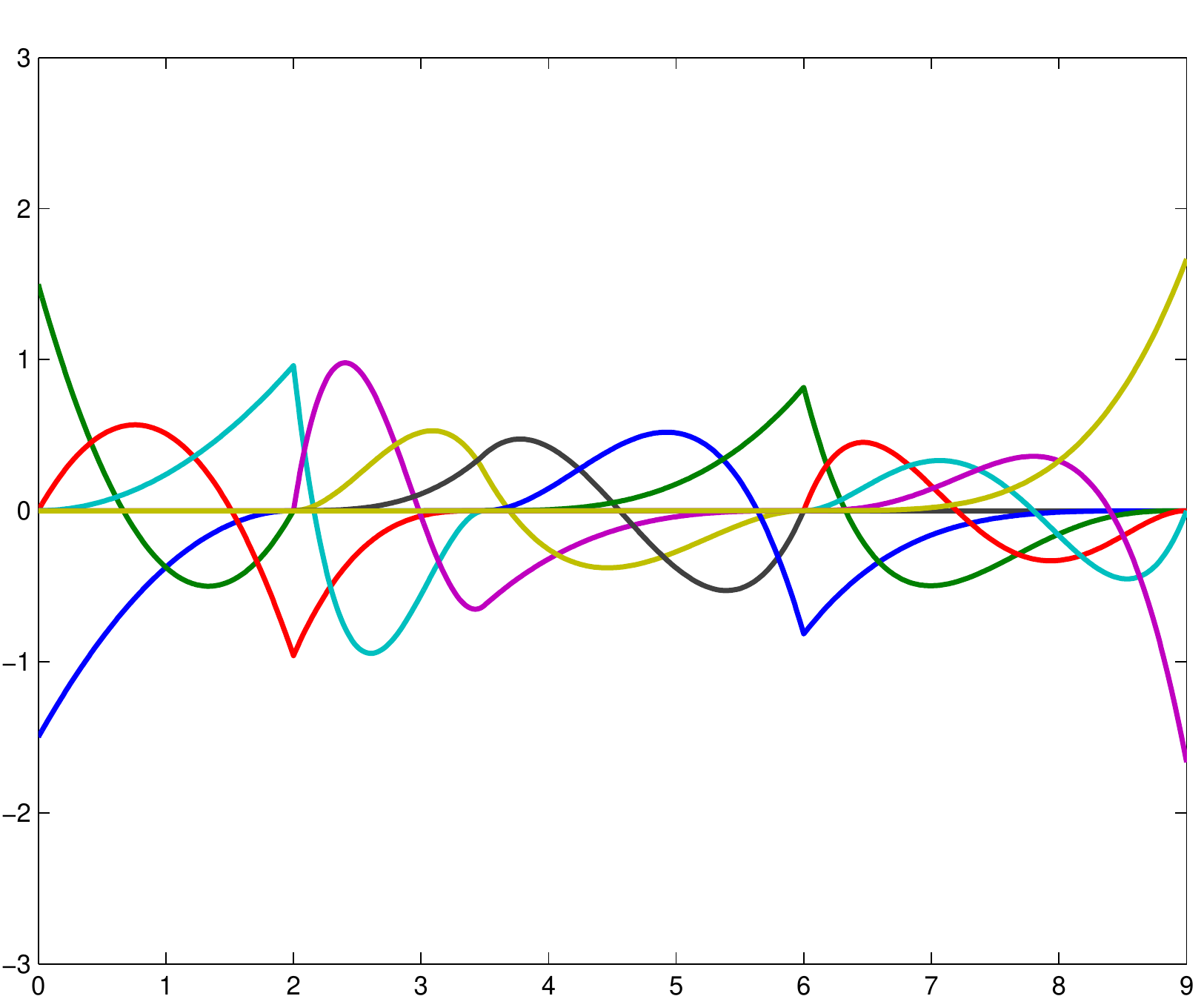}} \hspace*{0.1cm}
\subfigure[$\contin=2$]{\includegraphics[width=4.4cm]{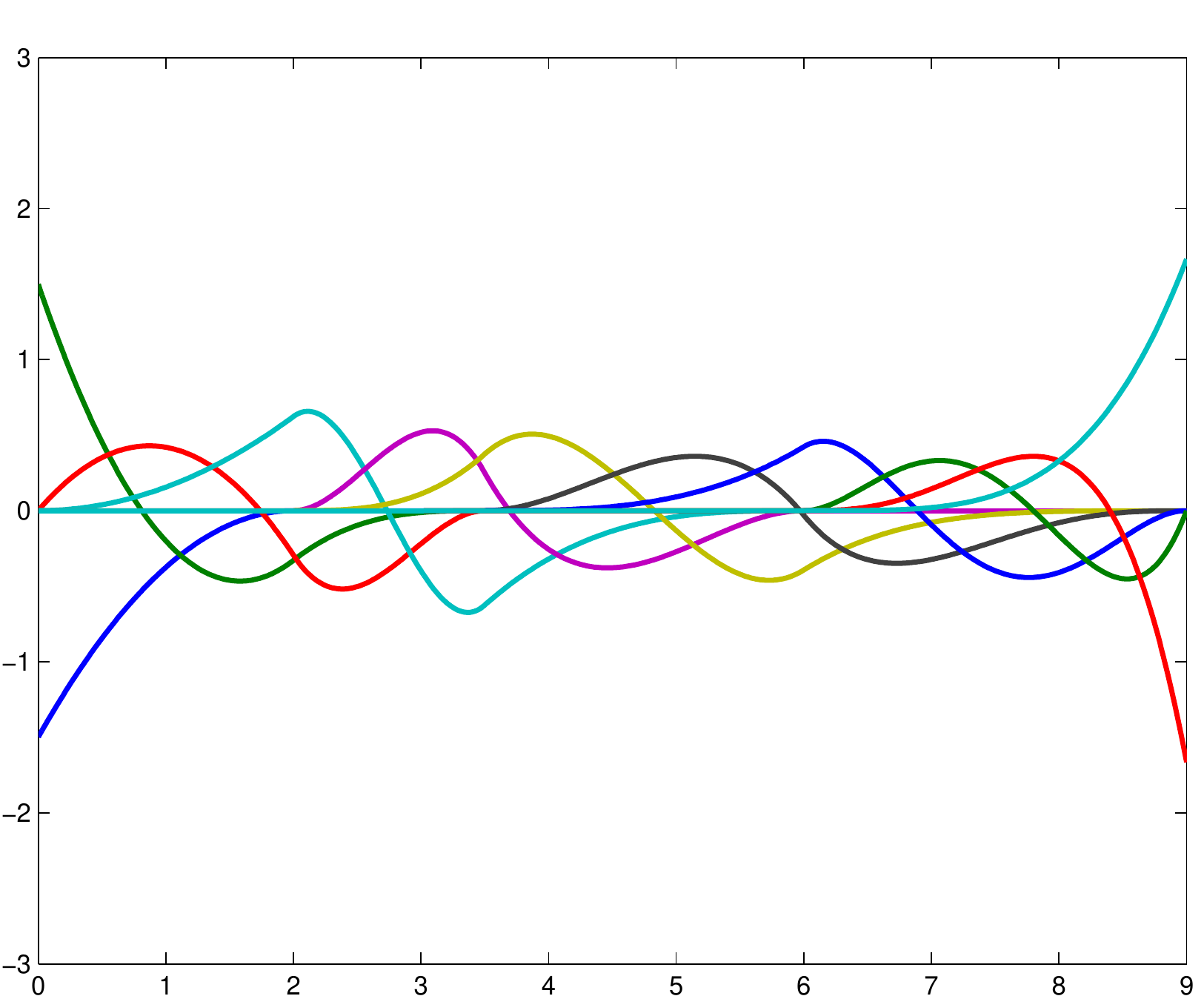}} \vspace*{-0.2cm}
\caption{First derivative of MDB-splines depicted in Figure~\ref{fig:ex-basis} }
\label{fig:ex-basis-der}
\end{figure}

\begin{example}\label{ex:ex-basis}
  Consider the following three open knot vectors:
  \begin{equation}\label{eq:ex1-knots}
  \begin{gathered}
    \seg{1}{U} = [0,\; 0,\; 0,\; 0,\; 2,\; 2,\; 2,\; 2], \\
    \seg{2}{U} = [0,\; 0,\; 0,\; 0,\; 0,\; {3}/{2},\; {3}/{2},\; 4,\; 4,\; 4,\; 4,\; 4], \\
    \seg{3}{U} = [0,\; 0,\; 0,\; 0,\; 0,\; 0,\; 3,\; 3,\; 3,\; 3,\; 3,\; 3].
  \end{gathered}
  \end{equation}
  They give rise to three B-spline bases of degree $\seg{1}{\pdeg}=3$, $\seg{2}{\pdeg}=4$, and $\seg{3}{\pdeg}=5$, respectively. Starting from these local B-spline bases,
  we create the MDB-spline basis associated with the space $\SPL{S}{\contin}{U}$ where $\mbf{U}=(\seg{1}{U},\seg{2}{U},\seg{3}{U})$ and $\mbf{\contin}=(\contin,\contin)$ for different choices of $\contin$. The resulting MDB-splines for $\contin\in\{0,1,2\}$ are shown in Figure~\ref{fig:ex-basis}, and their first derivative in Figure~\ref{fig:ex-basis-der}. 
  The smoothness $C^\kappa$ at each segment join $\xi\in\{2,6\}$ can be visually verified through the derivative plots. The dimension of the spline space equals $15-2\contin$, in accordance with the formula in \eqref{eq:dimension-md}.
\end{example}

\begin{figure}[t!]
\subfigure[Periodic MDB-splines]{\includegraphics[width=4.4cm]{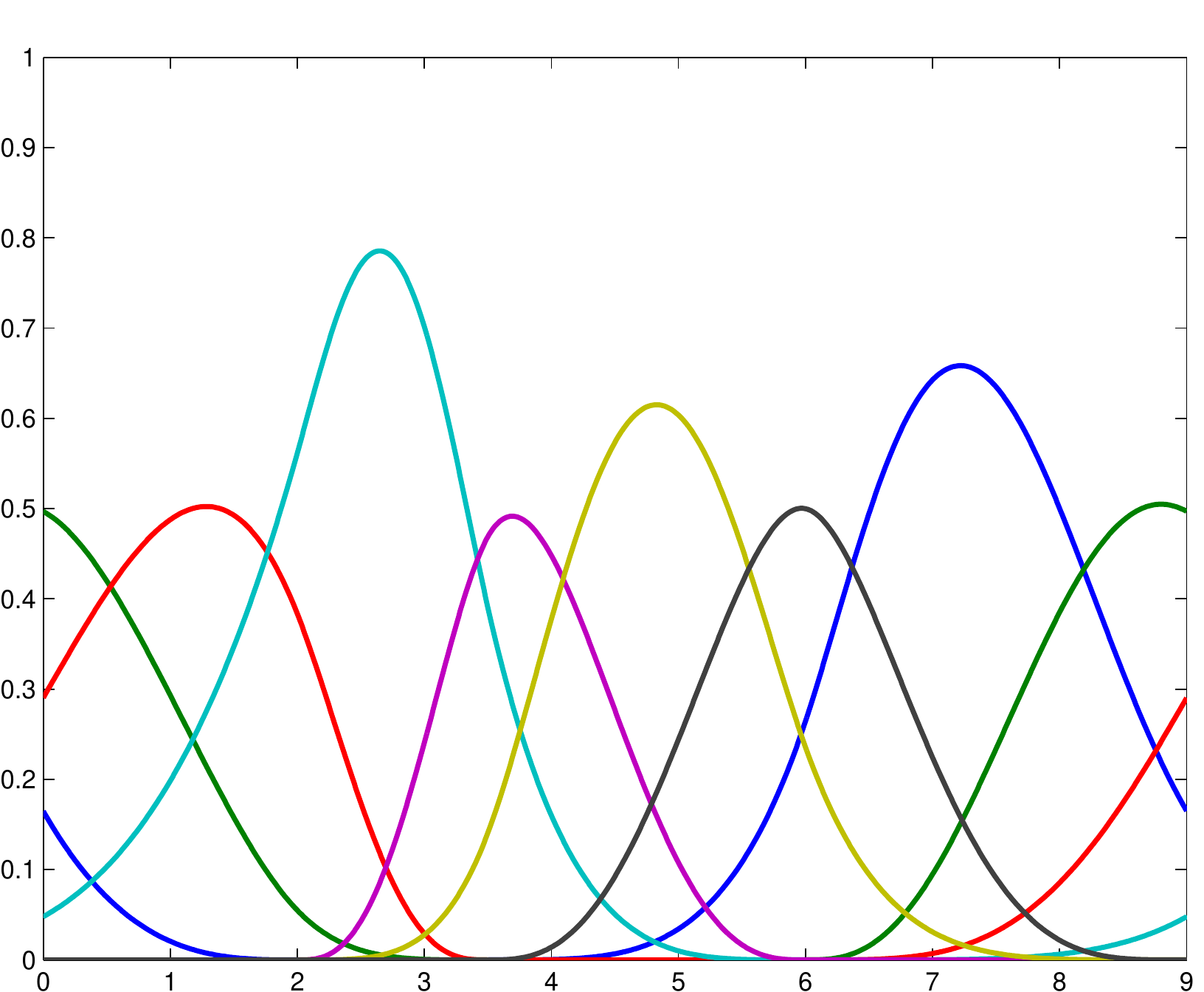}} \hspace*{0.1cm}
\subfigure[First derivative of MDB-splines]{\includegraphics[width=4.4cm]{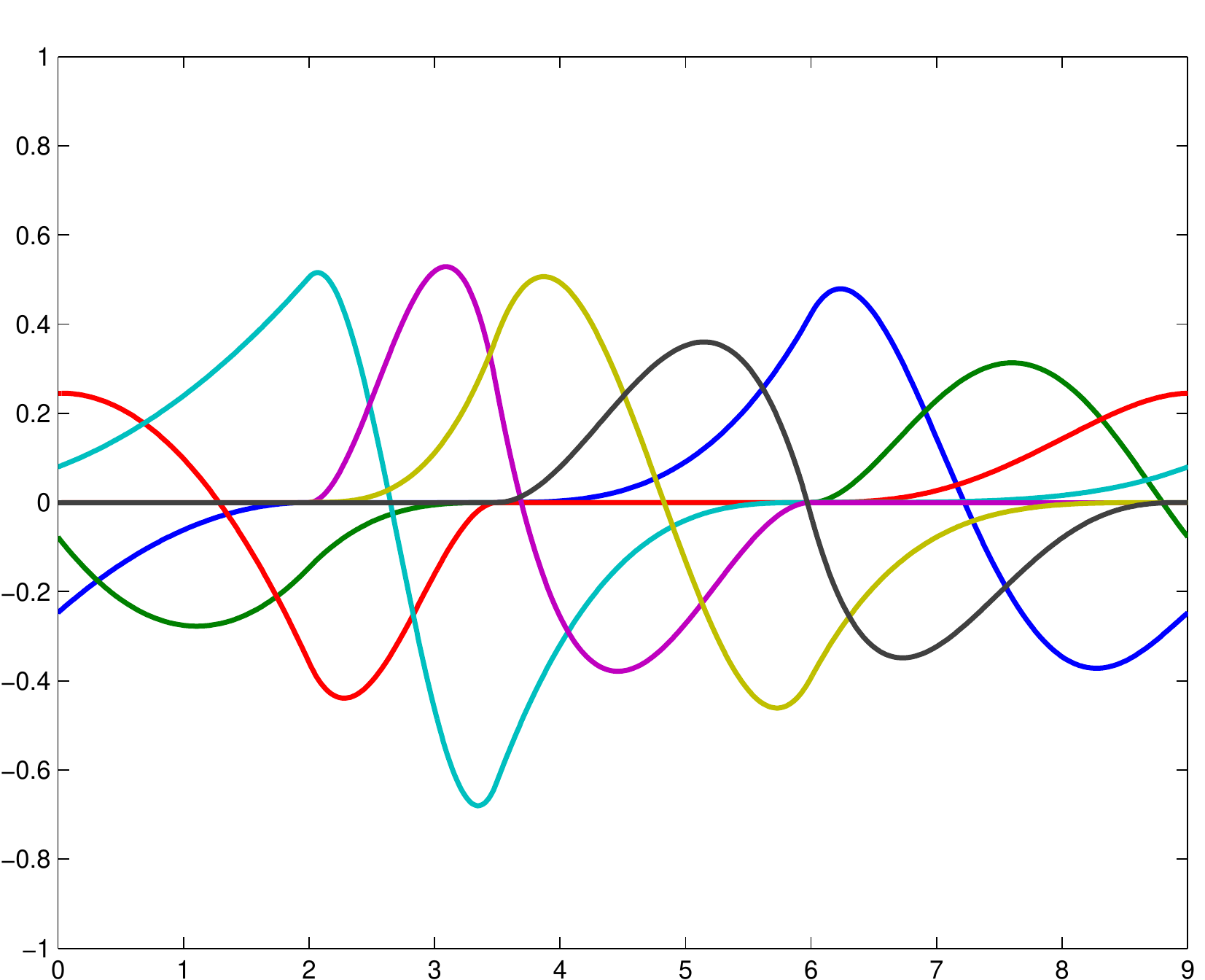}}  \vspace*{-0.2cm}
\caption{Periodic MDB-spline basis of degrees $(3,4,5)$ defined over the knot vector configuration $\mbf{U}=(\seg{1}{U},\seg{2}{U},\seg{3}{U})$ specified in \eqref{eq:ex1-knots} and continuity orders $\mbf{\contin}=(2,2)$ and $\per{\contin}=3$}
\label{fig:ex-basis-per}
\bigskip
\subfigure[$C^1$ spline function]{\includegraphics[width=4.4cm]{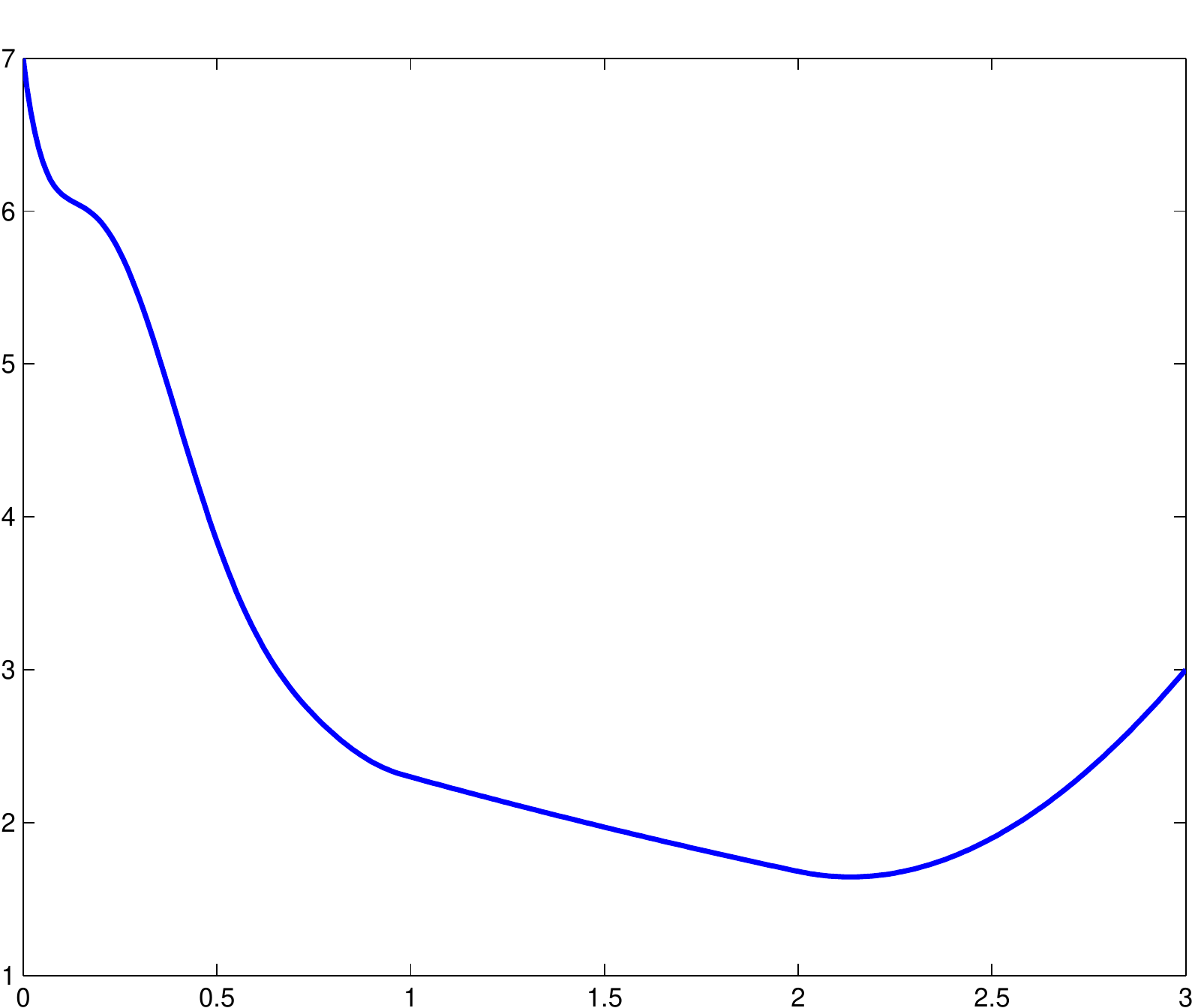}} \hspace*{0.1cm}
\subfigure[MDB-splines of degrees $(7,2,3)$]{\includegraphics[width=4.4cm]{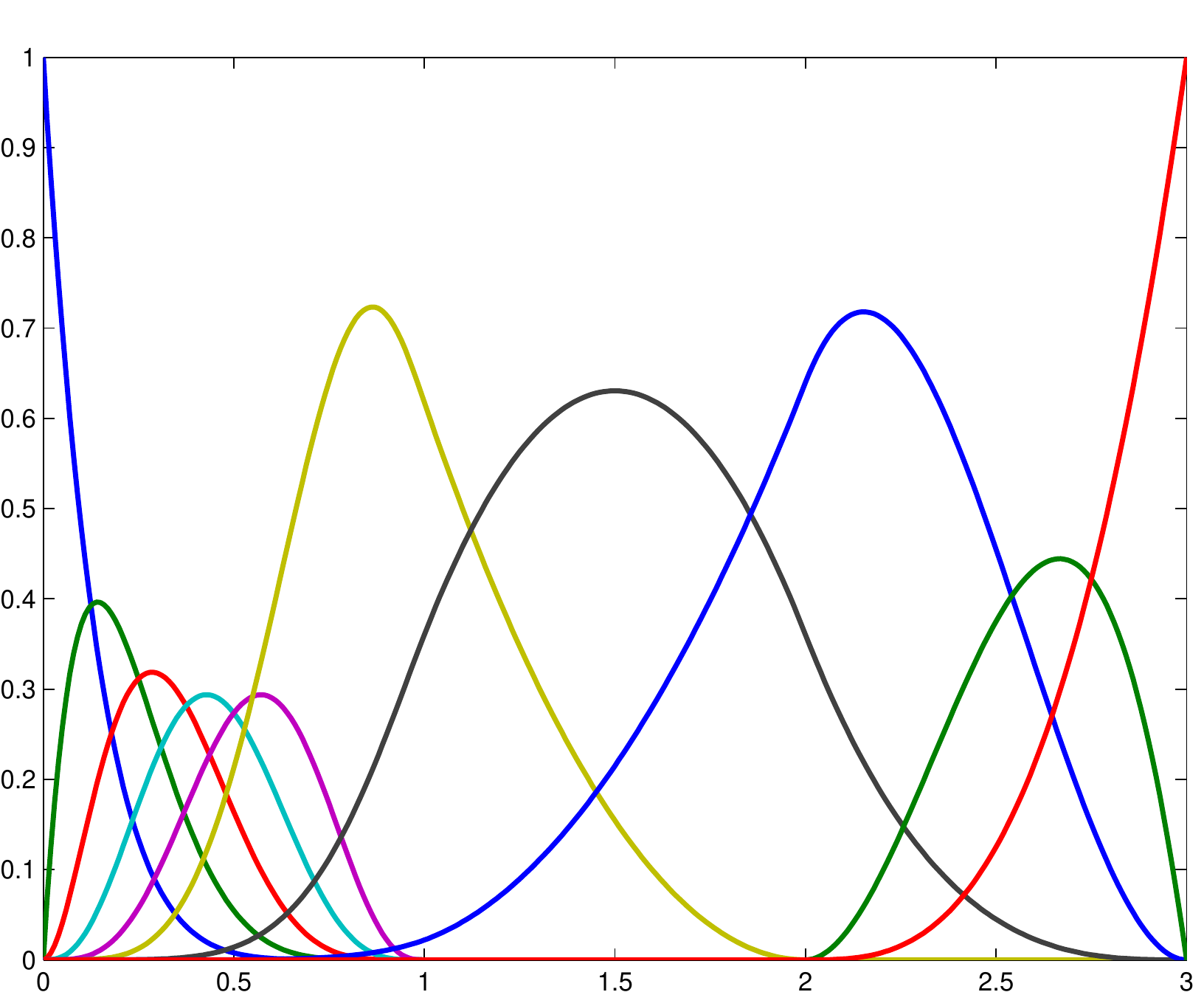}} \hspace*{0.1cm}
\subfigure[B-splines of degree $7$]{\includegraphics[width=4.4cm]{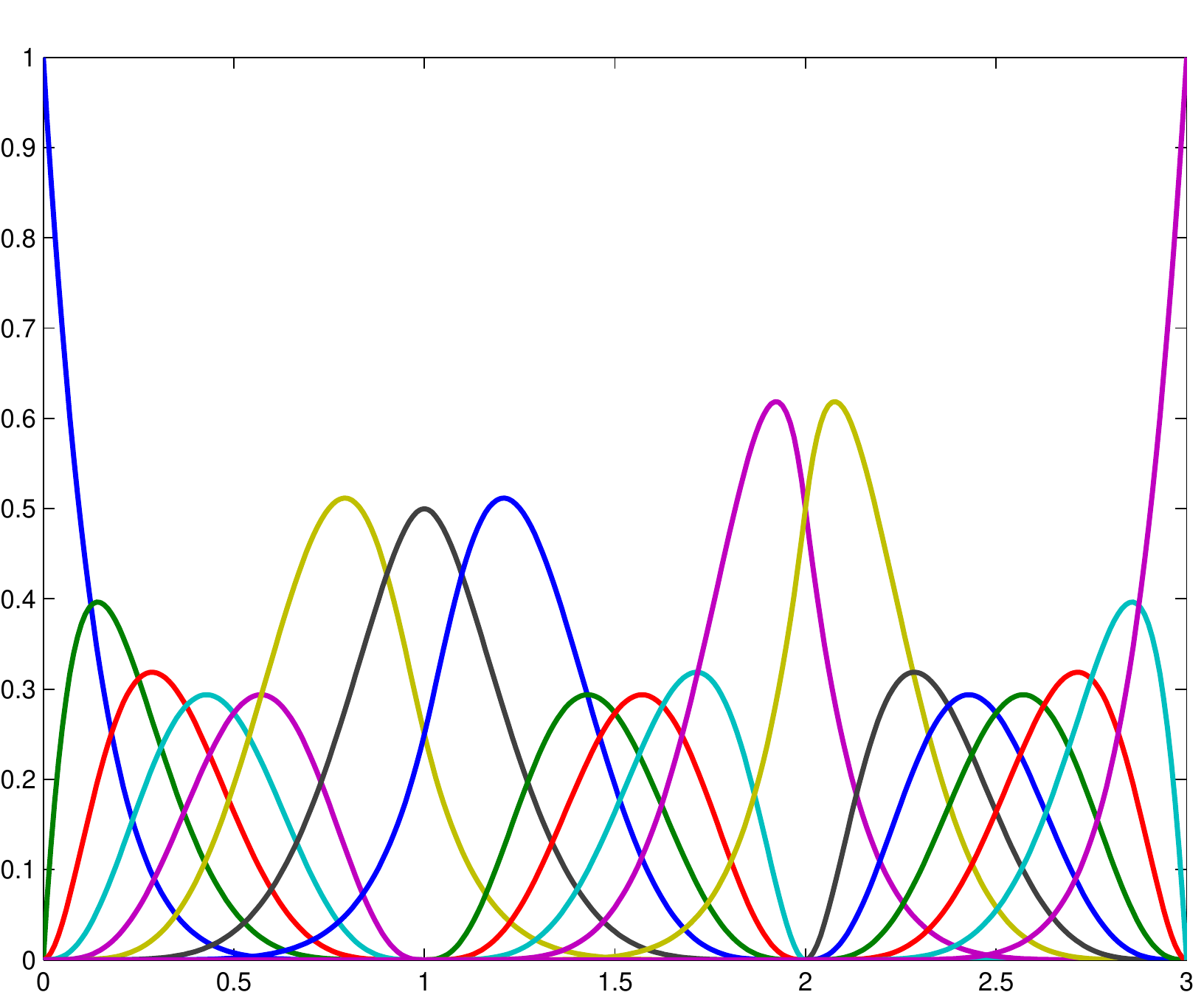}} \vspace*{-0.2cm}
\caption{A given $C^1$ spline function represented in terms of the MDB-spline basis over the knot vector configuration $\mbf{U}=(\seg{1}{U},\seg{2}{U},\seg{3}{U})$ specified in \eqref{eq:ex2-knots} and continuity orders $\mbf{\contin}=(2,1)$, together with the alternative B-spline basis over the knot vector specified in \eqref{eq:ex2-knots-full}}
\label{fig:ex-conversion}
\end{figure}

\begin{example}\label{ex:ex-basis-per}
  Consider again the three open knot vectors in \eqref{eq:ex1-knots}.
  We now create a periodic MDB-spline basis associated with the space $\SPL{S}{\contin}{U}$ where $\mbf{U}=(\seg{1}{U},\seg{2}{U},\seg{3}{U})$ and $\mbf{\contin}=(2,2)$; the periodic continuity order is chosen to be $\per{\contin}=3$. The resulting MDB-splines and their first derivative are shown in Figure~\ref{fig:ex-basis-per}.
  Note that the three central MDB-splines in Figure~\ref{fig:ex-basis-per}(a) coincide with the ones in Figure~\ref{fig:ex-basis}(c), illustrating the locality of the basis construction.
\end{example}

\begin{example}\label{ex:ex-conversion}
  Consider the following three open knot vectors:
  \begin{equation}\label{eq:ex2-knots}
  \begin{gathered}
    \seg{1}{U} = [0,\; 0,\; 0,\; 0,\; 0,\; 0,\; 0,\; 0,\; 1,\; 1,\; 1,\; 1,\; 1,\; 1,\; 1,\; 1], \\
    \seg{2}{U} = [0,\; 0,\; 0,\; 1,\; 1,\; 1], \quad
    \seg{3}{U} =  [0,\; 0,\; 0,\; 0,\; 1,\; 1,\; 1,\; 1].
  \end{gathered}
  \end{equation}
  They allow for constructing three B-spline bases of degree $\seg{1}{\pdeg}=7$, $\seg{2}{\pdeg}=2$, and $\seg{3}{\pdeg}=3$, respectively. Then, we focus on the multi-degree spline space $\SPL{S}{\contin}{U}$ with $\mbf{U}= (\seg{1}{U},\seg{2}{U},\seg{3}{U})$ and $\mbf{\contin}=(2,1)$, which has dimension $10$. Each spline in this space can be represented in terms of standard B-splines of degree $\pdeg=7$ (i.e., the highest MDB-spline degree) defined on the knot vector
  \begin{equation}\label{eq:ex2-knots-full}
    U = [0,\; 0,\; 0,\; 0,\; 0,\; 0,\; 0,\; 0,\; 1,\; 1,\; 1,\; 1,\; 1,\; 2,\; 2,\; 2,\; 2,\; 2,\; 2,\; 3,\; 3,\; 3,\; 3,\; 3,\; 3,\; 3,\; 3].
  \end{equation}
  The B-spline space has a much larger dimension, namely $19$. Conversion from a spline $s\in\SPL{S}{\contin}{U}$ in MDB-spline representation (of degrees $7,2,3$) to a standard B-spline representation (of degree $7$) can be done by local degree raising. 
  As an example, referring to Figure~\ref{fig:ex-conversion}, the $C^1$ spline function on the left can be represented in terms of the MDB-splines (middle) using the coefficient vector
  \begin{equation*}
   \mbf{s}_{\mathrm{MDB}} = [7,\; 4,\; 10,\; 1,\; 4,\; 2.5,\; 2,\; 1.5,\; 2,\; 3 ],
  \end{equation*}
  and in terms of the B-splines (right) using the coefficient vector
  \begin{align*}
       \mbf{s}_{\mathrm{B}} = [7,\; 4,\; 10,\; 1,\; 4,\; 2.5,\; 2.2941,\; 2.1029,\; 2.0110,\; 1.9228,\; 1.8382,\; & \\
        1.7574,\; 1.6029,\; 1.6229,\; 1.7349,\; 1.9337,\; 2.2143,\; 2.5714,\; 3]. &
  \end{align*}
  This indicates that the multi-degree structure could be a powerful tool for data compression by locally tuning the degree according to the data instead of taking a globally uniform degree.
\end{example}

\section{Conclusions} \label{sec:conclusion}

In the paper we have focused on the algorithmic aspects of multi-degree splines.
These are a natural extension of standard splines, consisting of polynomial pieces that can have different degrees. 
The flexibility of choosing the degree of each piece (or a cluster of pieces) independently of the others
allows us to simplify the global description of very complicated functions/profiles by requiring locally lower degrees, and hence less degrees-of-freedom.

Multi-degree splines admit a basis (the so-called MDB-spline basis) that can be regarded as a natural generalization of the standard B-spline basis. MDB-splines inherit many important properties, such as local support, non-negative partition of unity, and smoothness at segment joins. 
We have detailed a simple procedure to construct an extraction operator that represents all MDB-splines as linear combinations of local B-splines of different degrees. Such extraction operator enables the use of existing efficient algorithms for B-spline evaluations and refinements in the context of multi-degree splines.
A small \Matlab toolbox has been provided to illustrate the computation and use of MDB-splines (it is available through the CALGO library).

The framework for smoothly joining different (B-spline) spaces is quite abstract and can be applied in a more general  setting.
The main algorithm requires only information on the derivatives of the underlying basis functions at the joins (see Figure~\ref{alg:extraction}). Therefore, it could be easily adapted towards other basis functions as long as they have similar structural properties to B-splines. 
For example, in \cite{toshniwal2017polar} the technique was employed to smoothly join rational NURBS spaces. However, we recall that the extension to NURBS has to be applied with care because non-negativity is lost in general beyond the $C^1$ case.
The use of Tchebycheffian B-spline spaces (generated by different ECT-systems) instead of B-spline spaces could be another interesting extension, in the spirit of \cite{buchwald2003}.

Finally, we believe that the presented algorithmic approach may pave the path for multi-degree splines from a rather theoretical to a truly practical tool. Application fields of interest include CAGD, IGA, data fitting and compression. The potential of MDB-splines in CAGD is well explained and illustrated in \cite{beccari2017}. In particular, control points can be defined and manipulated similarly to standard B-splines. However, it still has to be investigated how the local degrees can be specified in a geometrically intuitive way to be fully mature for practical CAGD software. 
In \cite{toshniwal2017polar} it was shown that the MDB-spline approach through an extraction operator is very convenient in IGA since it enables the conversion from a non-smooth B-spline multi-patch representation to a smooth MDB-spline representation. The higher smoothness can then be beneficially employed for solving high-order differential problems.
Because MDB-splines allow for local degree raising, they also offer new possibilities in IGA for adaptive local $p$- and $k$-refinement. At the moment, such refinements are very promising techniques in IGA but can only be done globally using the standard B-spline framework 
\cite{cottrell2009}; 
MDB-splines offer here a possibly substantial computation saving thanks to the local degree flexibility. Because of the implementation in terms of an extraction operator, MDB-splines can be quite easily integrated in most of the existing IGA and FEA software libraries without the need for modifying their kernel assembly or data structures (especially when standard B-splines or Bernstein polynomials are already supported); we refer the reader to \cite{borden2011isogeometric} for more details. 
Last but not least, MDB-splines have the ability for more efficient data compression by adaptively choosing the degree according to the local difficulties of the given data. An optimal choice between the local mesh size and local degree is an interesting topic of further research.

\begin{acks}

This work was partially supported by the Mission Sustainability Programme of the University of Rome Tor Vergata through the project IDEAS (CUP E81I18000060005) and by the MIUR Excellence Department Project awarded to the Department of Mathematics, University of Rome Tor Vergata (CUP E83C18000100006). The author is a member of Gruppo Nazionale per il Calcolo Scientifico -- Istituto Nazionale di Alta Matematica. 

\end{acks}


\bibliographystyle{ACM-Reference-Format}
\bibliography{mdb_splines}

\end{document}